\documentclass[12pt]{article}
\usepackage[cp1251]{inputenc}
\usepackage[russian]{babel}
\usepackage{graphicx}
\usepackage{amsmath,amssymb,amsopn,eufrak,  amsthm,amsfonts}
\begin{document}
\hyphenation{}
\newtheorem{theo}{Теорема}
\newtheorem{pr}{Предложение}
\newtheorem{cor}{Следствие}
\newtheorem{lemma}{Лемма}
\renewcommand{\labelenumi}{(\theenumi)}

\begin{center}Е.В. Журавлев, А.С. Кузьмина, Ю.Н. Мальцев

\vspace{1cm}

{\bf ОПИСАНИЕ МНОГООБРАЗИЙ КОЛЕЦ, В КОТОРЫХ КОНЕЧНЫЕ КОЛЬЦА ОДНОЗНАЧНО ОПРЕДЕЛЯЮТСЯ СВОИМИ ГРАФАМИ ДЕЛИТЕЛЕЙ НУЛЯ}

\vspace{1cm}
\end{center}

\sloppy

В данной работе рассматриваются ассоциативные кольца (не обязательно коммутативные и не обязательно имеющие единицу).

{\bf Определение. } {\it Графом делителей нуля кольца} $R$ называется граф, вершинами которого являются все ненулевые делители нуля кольца (односторонние и двусторонние), причем две различные вершины $x,y$ соединяются ребром тогда и только тогда, когда $xy=0$ или $yx=0$.

Обычно граф делителей нуля кольца $R$ обозначается через $\Gamma(R).$ Мы также будем использовать это обозначение.

Понятие графа делителей нуля было введено в работе \cite{Beck}. И. Бек ввел это понятие для {\it коммутативного} кольца и вершинами графа делителей нуля считал все элементы кольца. В статье \cite{AL} определение было изменено: в качестве вершин графа делителей нуля коммутативного кольца авторы этой работы рассматривали лишь ненулевые делители нуля. Затем понятие графа делителей нуля было распространено и на {\it некоммутативный} случай (см., например, \cite{AM-314}).

Нетрудно привести примеры неизоморфных колец, графы делителей нуля которых равны. Например, если $A$~--- счетномерная алгебра с нулевым умножением над полем~$\mathbb{Z}_p$, а $B$~--- счетномерная алгебра с нулевым умножением над полем~$\mathbb{Z}_q$, где $p,q$~--- это различные простые числа, то $\Gamma(A)\cong \Gamma(B),$ но $A\not\cong B.$ Другими словами, даже в многообразии $var\left\langle xy=0\right\rangle$ существуют примеры бесконечных неизоморфных колец, графы делителей нуля которых имеют одинаковое строение. В связи с этим интерес представляет такой вопрос: при каких условиях из равенства графов делителей нуля следует изоморфизм колец? Некоторые результаты, дающие ответ на этот вопрос для коммутативных колец, были получены в работе \cite{AM-274}. В настоящей работе данная проблема исследуется на языке многообразий, а именно: исследуются многообразия ассоциативных колец, в которых каждое конечное кольцо однозначно определяется своим графом делителей нуля. Другими словами, изучаются свойства многообразия колец $\mathfrak{M}$, для которого из равенства $\Gamma(R)=\Gamma(S)$ для конечных  колец $R,S\in \mathfrak{M},$ следует, что $R\cong S.$ Ранее такие многообразия исследовались в работах~\cite{semr, aejm2012}. Однако полного описания получено не было. В настоящей же работе многообразия, в которых все конечные кольца однозначно определяются своими графами делителей нуля, полностью описаны.

Введем обозначения и понятия, используемые в настоящей работе. 

{\em Полным $n$--вершинным графом} $K_n$ называется граф (без петель и кратных ребер), все $n$ вершин которого смежны между собой. 

Пусть аддитивная группа кольца $R$ разлагается в прямую сумму своих ненулевых аддитивных подгрупп $A_i,$ где $i=1,\ldots,n$ и $n\geq 2,$ т.е. $R=A_1\stackrel{.}{+}\ldots \stackrel{.}{+} A_n.$ Если все подгруппы $A_i$ являются двусторонними идеалами кольца $R,$ то кольцо $R$ называется {\it разложимым} ( в обозначении ${R=A_1 \oplus \ldots \oplus A_n}$). 

Порядок конечного кольца $R$ мы будем обозначать через $\left|R\right|.$ Для любого элемента $a\in R$, где $R$ -- произвольное кольцо, будем использовать следующее обозначение: $ann(a)=\{x\in R; xa=ax=0\}.$ Для любых элементов $x,y$ кольца $R$ положим $[x,y]=xy-yx$ и $x\circ y=xy+yx.$  Через $\mathbb{Z}_n$ мы будем обозначать кольцо классов вычетов по модулю $n$. Для простого числа $p$ будем полагать, что $N_{0,p}=\left\langle a\right\rangle, pa=0, a^2=0.$

Пусть $\mathbb{Z}\left\langle X\right\rangle=\mathbb{Z}\left\langle x_1,x_2,\ldots\right\rangle$ -- свободное ассоциативное кольцо от счетного числа переменных $X=\{x_1,x_2,\ldots\}$ и $f(x_1,\ldots,x_d)\in \mathbb{Z}\left\langle X\right\rangle$. Многочлен  $f(x_1,\ldots,x_d)$ {\it существенно зависит }от $x_1,x_2,\ldots, x_d$, если $f(0,x_2,\ldots,x_d)=\ldots=f(x_1,\ldots,x_{d-1},0)=0$. Минимальная из степеней одночленов, входящих в запись  $f(x_1,\ldots,x_d)$ с ненулевым коэффициентом, называется {\it нижней степенью многочлена}  $f(x_1,\ldots,x_d)$.

Пусть $\mathfrak{M}$ -- многообразие колец. Через $T(\mathfrak{M})$ будем обозначать множество всех многочленов из $\mathbb{Z}\left\langle X\right\rangle$, являющихся тождествами на всех кольцах из $\mathfrak{M}$. Назовем множество $T(\mathfrak{M})$ {\it идеалом тождеств} многообразия $\mathfrak{M}$. Если идеал тождеств $T(\mathfrak{M})$ порождается (как вполне характеристический идеал) многочленами $f_i,i\in I,$ то будем использовать следующее обозначение: $T(\mathfrak{M})=\{f_i~~|~~i\in I\}^T$. Через $\mathfrak{M}\vee \mathfrak{N}$ обозначается объединение многообразий  $\mathfrak{M}$ и $\mathfrak{N}$. Нетрудно заметить, что $T\left(\mathfrak{M}\vee \mathfrak{N}\right)=T\left(\mathfrak{M}\right)\cap T \left(\mathfrak{N}\right).$ 

Пусть $\mathfrak{M}_{1,p}=var\left\langle xyz=0, x^2=0, px=0\right\rangle$ и $\mathfrak{M}_{2,p}=var\left\langle xyz=0, [x,y]=0, px=0\right\rangle$, где $p$~--- произвольное простое число. Заметим, что $\mathfrak{M}_{1,2}\subseteq\mathfrak{M}_{2,2},$ а при нечетном $p$ мы имеем $\mathfrak{M}_{1,p}\cap\mathfrak{M}_{2,p}=var\left\langle xy=0, px=0\right\rangle$. Далее, пусть $F_i$~--- приведенно свободное кольцо с шестью порождающими $\{x_1, \ldots, x_6\}$ многообразия $\mathfrak{M}_{i,p},$ $i=1,2.$ Рассмотрим кольца $$A_{i,p}=F_i/\left\langle x_3x_4-x_1x_2\right\rangle,  B_{i,p}=F_i/\left\langle x_5x_6-x_1x_2-x_3x_4\right\rangle,$$ где  $\left\langle a,b,\ldots \right\rangle$~--- это идеал кольца $F_i,$ порожденный элементами $\{a,b, \ldots\},$ $i=1,2.$ Наша ближайшая цель~--- доказать, что $\Gamma(A_{1,p})\cong \Gamma(B_{1,p})$ для любого простого числа~$p$ и $\Gamma(A_{2,p})\cong \Gamma(B_{2,p})$ при нечетном простом~$p.$

\vspace{0.5cm}

По теореме Тарского любое ненулевое многообразие колец содержит одно из  минимальных многообразий:  $var~\mathbb{Z}_p$ или $var~N_{0,p},$ где $p$~-- некоторое простое число~\cite{Tarski}. Оказывается, что в минимальных многообразиях $var~\mathbb{Z}_p$ и $var~N_{0,p},$ где $p$~-- любое простое число, все конечные кольца однозначно определяются своими графами делителей нуля так же, как и в многообразии $var~N_{0,p_1}\vee \ldots \vee var~N_{0,p_s}$ \cite{semr}. В многообразии ${var~N_{0,p_1}\vee \ldots \vee var~N_{0,p_s}\vee var~\mathbb{Z}_p ,}$ где $p_1, \ldots, p_s$~-- попарно различные простые числа, $p$~-- любое простое число (возможно, совпадающее с одним из чисел~$p_i$), каждое конечное кольцо однозначно определяется своим графом делителей нуля тогда и только тогда, когда $(p_i,p)\neq(3,2)$ при~$i\leq s$ \cite{semr}. Наша цель --- показать, что любое многообразие, в котором все конечные кольца однозначно определяются своими графами делителей нуля, является подмногообразием многообразия вида ${var~N_{0,p_1}\vee \ldots \vee var~N_{0,p_s}\vee var~\mathbb{Z}_p},$ где $p_1, \ldots, p_s$~-- попарно различные простые числа, $p$~-- любое простое число и ${(p_i,p)\neq(3,2)}$ при всех~$i\leq s$. Для этого нам понадобятся некоторые вспомогательные утверждения. 

Покажем сначала, что $\Gamma(A_{1,p})\cong \Gamma(B_{1,p})$, в то время, как $A_{1,p}\not\cong B_{1,p}$ для любого простого числа~$p$. 

\textbf{Замечание. } Порядок алгебр $A_{i,p}, B_{i,p},$ $i=1,2,$ равен $p^{14}$. Возникает вопрос  о  существовании в указанных многообразиях примеров неизоморфных конечных колец небольшого порядка с одинаковыми графами делителей нуля. Однако нами было доказано, что в многообразии~$\mathfrak{M}_{1,2}$ все конечные кольца порядка $\leq 64$ однозначно определяются своими графами делителей нуля (фактически были полностью описаны все конечные кольца в многообразии~$\mathfrak{M}_{1,2}$, порядок которых не превышает~$64$).

\vspace{0.5cm}

\begin{lemma}\label{l11}Множество $C_1=\{\overline{x}_i\overline{x}_j;~ (i,j)\neq(3,4), 1\leq i<j\leq 6\}$ является базисом алгебры $A_{1,p}^2$, где $p$~--- простое число.)\end{lemma}
\begin{proof} Множество~$C_1$ является системой образующих векторного пространства~$A_{1,p}^2.$ Докажем, что это множество линейно независимо. Если множество~$C_1$ линейно зависимо, то  существуют элементы $\alpha,  \delta_{ij}\in\mathbb{Z}_p,$ не все равные нулю, такие, что в алгебре~$F_1$ справедливо  равенство \begin{equation}\sum\limits_{\tiny\begin{array}[c]{c}  i<j \\ (i,j)\neq (3,4)\end{array}}\delta_{ij}x_ix_j=\alpha(x_3x_4-x_1x_2).\end{equation} Положим  в этом равенстве ${x_3=x_4=x_5=x_6=0.}$ Тогда $\delta_{12}=-\alpha,$ другими словами, $\sum\limits_{\tiny\begin{array}[c]{c}  i<j \\ (i,j)\neq (3,4)\\ (i,j)\neq (1,2) \end{array}}\delta_{ij}x_ix_j=\alpha x_3x_4.$ Положим ${x_1=x_2=x_5=x_6=0.}$ Получим, что $\alpha=0$ и $\delta_{ij}=0$ для всех $i,j,$ таких, что $i<j $ и $(i,j)\neq (3,4), (1,2)$. Противоречие доказывает лемму. \end{proof}

\begin{pr}\label{pr1} Пусть $a=\sum\limits^{6}_{i=1}\alpha_i\overline{x}_i+u,$ $b=\sum\limits^{6}_{i=1}\beta_i\overline{x}_i+v$~--- произвольные элементы из кольца $A_{1,p},$ причем $u,v\in A_{1,p}^2.$ Если $a,b\notin A_{1,p}^2,$ то $ab=\overline{0}$ тогда и только тогда, когда $(\beta_1, \ldots, \beta_6)=\lambda(\alpha_1, \ldots, \alpha_6)$ для некоторого ненулевого элемента $\lambda\in\mathbb{Z}_p.$ \end{pr}
\begin{proof} 
Пусть $(\beta_1, \ldots, \beta_6)=\lambda(\alpha_1, \ldots, \alpha_6)$, где $0\neq\lambda\in\mathbb{Z}_p.$ Тогда $a=c+u,$ $b=\lambda c+v$ и $ab=\lambda c^2=\overline{0},$ т.к. $x^2=0$~--- тождество в алгебре $A_{1,p}$.

Докажем обратное утверждение. Пусть $ab=\overline{0}.$ Тогда $$ab=\left(\left|\begin{array}[c]{cc}  \alpha_1 &  \alpha_2 \\ \beta_1 & \beta_2 \end{array}\right|+\left|\begin{array}[c]{cc}  \alpha_3 &  \alpha_4 \\ \beta_3 & \beta_4 \end{array}\right|\right)\overline{x}_1\overline{x}_2+\sum\limits_{\tiny\begin{array}[c]{c}  i<j \\ (i,j)\neq(1,2) \\ (i,j)\neq (3,4) \end{array}}\left|\begin{array}[c]{cc}  \alpha_i &  \alpha_j \\ \beta_i & \beta_j \end{array}\right|\overline{x}_i\overline{x}_j=\overline{0}.$$
По лемме~\ref{l11} множество $C_1=\{\overline{x}_i\overline{x}_j;~ (i,j)\neq(3,4), i<j\}$~--- базис векторного пространства $A_{1,p}^2.$ Поэтому получаем, что $$\left|\begin{array}[c]{cc}  \alpha_1 &  \alpha_2 \\ \beta_1 & \beta_2 \end{array}\right|+\left|\begin{array}[c]{cc}  \alpha_3 &  \alpha_4 \\ \beta_3 & \beta_4 \end{array}\right|=0 \mbox{~и~} \left|\begin{array}[c]{cc}  \alpha_i &  \alpha_j \\ \beta_i & \beta_j \end{array}\right|=0,$$ где $i<j,  (i,j)\neq(1,2)$ и $(i,j)\neq (3,4)$.

Рассмотрим следующие случаи.

\textbf{Случай 1. } Пусть $\alpha_1\neq 0.$

Поскольку $$\left|\begin{array}[c]{cc}  \alpha_1 &  \alpha_3 \\ \beta_1 & \beta_3 \end{array}\right|=0, \left|\begin{array}[c]{cc}  \alpha_1 &  \alpha_4 \\ \beta_1 & \beta_4 \end{array}\right|=0,$$ то вторые строки этих определителей линейно выражаются через первые, т.е. $\beta_1=\lambda \alpha_1,$ $\beta_3=\lambda \alpha_3$ и $\beta_1=\mu \alpha_1,$ $\beta_4=\mu\alpha_4$ для некоторых элементов $\lambda,\mu\in\mathbb{Z}_p.$ Поскольку $\alpha_1\neq 0$ и $(\lambda-\mu)\alpha_1=\beta_1-\beta_1=0,$ то $\lambda=\mu$ и $\beta_4=\lambda\alpha_4$. Аналогично доказывается, что $\beta_5=\lambda\alpha_5$ и $\beta_6=\lambda\alpha_6$.

Далее, из равенств $\left|\begin{array}[c]{cc}  \alpha_1 &  \alpha_2 \\ \beta_1 & \beta_2 \end{array}\right|+\left|\begin{array}[c]{cc}  \alpha_3 &  \alpha_4 \\ \beta_3 & \beta_4 \end{array}\right|=0$, $\beta_3=\lambda\alpha_3$ и $\beta_4=\lambda\alpha_4$ следует, что $\left|\begin{array}[c]{cc}  \alpha_1 &  \alpha_2 \\ \beta_1 & \beta_2 \end{array}\right|=0$ и, значит, $\beta_2=\lambda\alpha_2.$ Таким образом, $(\beta_1, \ldots, \beta_6)=\lambda(\alpha_1, \ldots, \alpha_6)$.

\textbf{Случай 2. } Пусть $\alpha_1= 0.$

Если хотя бы один из элементов $\alpha_2, \alpha_3, \alpha_4$ не равен нулю, то рассуждаем так же, как при рассмотрении случая~1. Поэтому можем полагать, что $\alpha_1=\alpha_2=\alpha_3=\alpha_4=0$ и, например, $\alpha_5\neq 0.$ (Заметим, что поскольку $a\notin A_{1,p}^2$, то один из элементов $\alpha_5$ или $\alpha_6$ отличен от нуля.) Аналогично мы можем считать, что $\beta_1=\beta_2=\beta_3=\beta_4=0.$ Следовательно, $$ab=\left|\begin{array}[c]{cc}  \alpha_5 &  \alpha_6 \\ \beta_5 & \beta_6 \end{array}\right|\overline{x}_5\overline{x}_6=\overline{0}, \mbox{~или~} \left|\begin{array}[c]{cc}  \alpha_5 &  \alpha_6 \\ \beta_5 & \beta_6 \end{array}\right|=0,$$ т.е. $\beta_5=\lambda \alpha_5, $ $\beta_6=\lambda \alpha_6$ для некоторого ненулевого элемента $\lambda\in\mathbb{Z}_p.$ Таким образом, $(\beta_1, \ldots, \beta_6)=\lambda(\alpha_1, \ldots, \alpha_6)$.\end{proof}

Из предложения~\ref{pr1} получаем, что для каждого элемента ${a=\sum\limits^{6}_{i=1}\alpha_i\overline{x}_i\in A_{1,p}\setminus A_{1,p}^2}$ множество $\{\alpha  a+u;~  0\neq \alpha\in \mathbb{Z}_p, u\in A_{1,p}^2\}$ образует полный подграф $\Gamma_a$ графа~$\Gamma(A_{1,p})$, причем каждая вершина подграфа $\Gamma_a$ смежна только с вершинами этого подграфа и элементами из множества~$A_{1,p}^2$. Покажем, что такое же строение имеет граф $\Gamma(B_{1,p})$. 

\begin{lemma}\label{l12}Множество $D_1=\{ \overline{x}_i\overline{x}_j;~ (i,j)\neq(5,6), 1\leq i<j\leq 6\}$ является базисом алгебры $B_{1,p}^2$ для любого простого числа $p$.\end{lemma}

\textit{Доказательство} аналогично доказательству леммы~\ref{l11}.

\begin{pr}\label{pr2}Пусть $a=\sum\limits^{6}_{i=1}\alpha_i\overline{x}_i+u,$ $b=\sum\limits^{6}_{i=1}\beta_i\overline{x}_i+v$~--- произвольные элементы из кольца $B_{1,p},$ где $u,v\in B_{1,p}^2.$ Если $a,b\notin B_{1,p}^2,$ то $ab=\overline{0}$ тогда и только тогда, когда $(\beta_1, \ldots, \beta_6)=\lambda(\alpha_1, \ldots, \alpha_6)$ для некоторого ненулевого элемента $\lambda\in\mathbb{Z}_p.$  \end{pr}
\begin{proof} 
Пусть $(\beta_1, \ldots, \beta_6)=\lambda(\alpha_1, \ldots, \alpha_6)$, где $0\neq\lambda\in\mathbb{Z}_p.$ Тогда $a=c+u,$ $b=\lambda c+v$ и $ab=\lambda c^2=\overline{0},$ т.к. $x^2=0$ для любого элемента $x\in B_{1,p}$.

Докажем обратное утверждение. Пусть $ab=\overline{0}.$ Тогда 
\begin{equation*}\begin{split}
ab=\left(\left|\begin{array}[c]{cc}  \alpha_1 &  \alpha_2 \\ \beta_1 & \beta_2 \end{array}\right|+\left|\begin{array}[c]{cc}  \alpha_5 &  \alpha_6 \\ \beta_5 & \beta_6 \end{array}\right|\right)\overline{x}_1\overline{x}_2+&\left(\left|\begin{array}[c]{cc}  \alpha_3 &  \alpha_4 \\ \beta_3 & \beta_4 \end{array}\right|+\left|\begin{array}[c]{cc}  \alpha_5 &  \alpha_6 \\ \beta_5 & \beta_6 \end{array}\right|\right)\overline{x}_3\overline{x}_4+ \\ +&\sum\limits_{\tiny\begin{array}[c]{c}  i<j \\ (i,j)\neq(1,2) \\ (i,j)\neq (3,4)\\ (i,j)\neq (5,6) \end{array}}\left|\begin{array}[c]{cc}  \alpha_i &  \alpha_j \\ \beta_i & \beta_j \end{array}\right|\overline{x}_i\overline{x}_j=0.\end{split}\end{equation*} По лемме~\ref{l12} множество $D_1=\{\overline{x}_i\overline{x}_j;~ (i,j)\neq(5,6), i<j\}$~--- базис векторного пространства $B_{1,p}^2.$ Поэтому  $$\left|\begin{array}[c]{cc}  \alpha_1 &  \alpha_2 \\ \beta_1 & \beta_2 \end{array}\right|+\left|\begin{array}[c]{cc}  \alpha_5 &  \alpha_6 \\ \beta_5 & \beta_6 \end{array}\right|=0, \left|\begin{array}[c]{cc}  \alpha_3 &  \alpha_4 \\ \beta_3 & \beta_4 \end{array}\right|+\left|\begin{array}[c]{cc}  \alpha_5 &  \alpha_6 \\ \beta_5 & \beta_6 \end{array}\right|=0 \mbox{~и~} \left|\begin{array}[c]{cc}  \alpha_i &  \alpha_j \\ \beta_i & \beta_j \end{array}\right|=0,$$ если $i<j,  (i,j)\neq(1,2), (i,j)\neq (3,4)$ и $(i,j)\neq (5,6)$.

Так же, как при доказательстве предложения~\ref{pr1}, рассмотрим два случая.

\textbf{Случай 1. } Пусть $\alpha_1\neq 0.$

Аналогично тому, как это было сделано в доказательстве предложения~\ref{pr1}, получаем, что $\beta_i=\lambda\alpha_i$ для любого  $i\neq2$. Поэтому из равенства ${\left|\begin{array}[c]{cc}  \alpha_1 &  \alpha_2 \\ \beta_1 & \beta_2 \end{array}\right|+\left|\begin{array}[c]{cc}  \alpha_5 &  \alpha_6 \\ \beta_5 & \beta_6 \end{array}\right|=0}$ вытекает равенство ${\left|\begin{array}[c]{cc}  \alpha_1 &  \alpha_2 \\ \beta_1 & \beta_2 \end{array}\right|=0}$. Отсюда следует, что $\beta_2=\lambda\alpha_2.$ Таким образом, при $\alpha_1\neq 0$ мы имеем, что $(\beta_1, \ldots, \beta_6)=\lambda(\alpha_1, \ldots, \alpha_6)$.

\textbf{Случай 2. } Пусть $\alpha_1= 0.$

Если хотя бы один из элементов $\alpha_2, \alpha_3, \alpha_4, \beta_1, \beta_2, \beta_3,\beta_4$ не равен нулю, то рассуждаем так же, как в случае~1. Поэтому можем считать, что $\alpha_i=\beta_j=0$ при $i,j\in\{1,2,3,4\}$. Значит, $$ab=(\alpha_5\overline{x}_5+\alpha_6\overline{x}_6)(\beta_5\overline{x}_5+\beta_6\overline{x}_6)=\left|\begin{array}[c]{cc}  \alpha_5 &  \alpha_6 \\ \beta_5 & \beta_6 \end{array}\right|(\overline{x}_1\overline{x}_2+\overline{x}_3\overline{x}_4)=\overline{0}.$$ Поскольку $a \notin B_{1,p}^2,$ то, например, $\alpha_5\neq 0.$ Таким образом, $\beta_5=\lambda\alpha_5,$ $\beta_6=\lambda\alpha_6$ для некоторого ненулевого элемента $\lambda\in\mathbb{Z}_p,$ т.е. $(\beta_1, \ldots, \beta_6)=\lambda(\alpha_1, \ldots, \alpha_6)$.
\end{proof}

Таким образом, справедливо

\begin{cor}\label{cor1}$\Gamma(A_{1,p})\cong\Gamma(B_{1,p})$ для любого простого числа~$p$.\end{cor}

\begin{pr}\label{pr3}Алгебра $A_{1,p}$ не изоморфна алгебре $B_{1,p}$  для любого простого числа~$p$.\end{pr}
\begin{proof} Предположим противное: существует изоморфизм $\varphi$ алгебры $B_{1,p}$ на алгебру $A_{1,p}$.
Введем такие обозначения:
$\overline{x_i}=x_i+\langle x_1x_2+x_3x_4-x_5x_6\rangle\in B_{1,p}$ и
${\overline{x}'_i=x_i+\langle x_1x_2-x_3x_4\rangle\in A_{1,p}}$, $i=\overline{1,6}$.
Заметим что $\varphi(B_{1,p}^2)\subseteq A_{1,p}^2$, а значит, отображение $\overline{\varphi}: B_{1,p}/B_{1,p}^2 \rightarrow A_{1,p}/A_{1,p}^2$, определенное по правилу
${\overline{\varphi}(a+B_{1,p}^2)=}$ ${=\varphi(a)+A_{1,p}^2}$, $a\in B_{1,p}$, также является изоморфизмом.
Отсюда следует, что существует невырожденная матрица $P=(p_{ij})_{6\times 6}$, $p_{ij}\in \mathbb{Z}_p,$ такая, что $$
\left\{
  \begin{array}{ll}
\varphi(\overline{x}_1)=p_{11}\overline{x}'_1+p_{21}\overline{x}'_2+p_{31}\overline{x}'_3+p_{41}\overline{x}'_4+p_{51}\overline{x}'_5+p_{61}\overline{x}'_6+b_1;\\
\varphi(\overline{x}_2)=p_{12}\overline{x}'_1+p_{22}\overline{x}'_2+p_{32}\overline{x}'_3+p_{42}\overline{x}'_4+p_{52}\overline{x}'_5+p_{62}\overline{x}'_6+b_2;\\
\varphi(\overline{x}_3)=p_{13}\overline{x}'_1+p_{23}\overline{x}'_2+p_{33}\overline{x}'_3+p_{43}\overline{x}'_4+p_{53}\overline{x}'_5+p_{63}\overline{x}'_6+b_3;\\
\varphi(\overline{x}_4)=p_{14}\overline{x}'_1+p_{24}\overline{x}'_2+p_{34}\overline{x}'_3+p_{44}\overline{x}'_4+p_{54}\overline{x}'_5+p_{64}\overline{x}'_6+b_4;\\
\varphi(\overline{x}_5)=p_{15}\overline{x}'_1+p_{25}\overline{x}'_2+p_{35}\overline{x}'_3+p_{45}\overline{x}'_4+p_{55}\overline{x}'_5+p_{65}\overline{x}'_6+b_5;\\
\varphi(\overline{x}_6)=p_{16}\overline{x}'_1+p_{26}\overline{x}'_2+p_{36}\overline{x}'_3+p_{46}\overline{x}'_4+p_{56}\overline{x}'_5+p_{66}\overline{x}'_6+b_6,
  \end{array}
\right.
$$
где $b_1$, $b_2$, $b_3$, $b_4$, $b_5$, $b_6$ -- некоторые элементы из $A_{1,p}^2$.

Учитывая свойства изоморфизма, получаем
\begin{gather*}
\varphi(\overline{x}_1\overline{x}_2+\overline{x}_3\overline{x}_4-\overline{x}_5\overline{x}_6)=\\
=(p_{11}\overline{x}'_1+p_{21}\overline{x}'_2+p_{31}\overline{x}'_3+p_{41}\overline{x}'_4+p_{51}\overline{x}'_5+p_{61}\overline{x}'_6)
(p_{12}\overline{x}'_1+p_{22}\overline{x}'_2+p_{32}\overline{x}'_3+p_{42}\overline{x}'_4+p_{52}\overline{x}'_5+p_{62}\overline{x}'_6)+\\
+(p_{13}\overline{x}'_1+p_{23}\overline{x}'_2+p_{33}\overline{x}'_3+p_{43}\overline{x}'_4+p_{53}\overline{x}'_5+p_{63}\overline{x}'_6)
(p_{14}\overline{x}'_1+p_{24}\overline{x}'_2+p_{34}\overline{x}'_3+p_{44}\overline{x}'_4+p_{54}\overline{x}'_5+p_{64}\overline{x}'_6)-\\
-(p_{15}\overline{x}'_1+p_{25}\overline{x}'_2+p_{35}\overline{x}'_3+p_{45}\overline{x}'_4+p_{55}\overline{x}'_5+p_{65}\overline{x}'_6)
(p_{16}\overline{x}'_1+p_{26}\overline{x}'_2+p_{36}\overline{x}'_3+p_{46}\overline{x}'_4+p_{56}\overline{x}'_5+p_{66}\overline{x}'_6)=\\
=(p_{11}p_{62}-p_{12}p_{61}+p_{13}p_{64}-p_{14}p_{63}-p_{15}p_{66}+p_{16}p_{65})\overline{x}'_1\overline{x}'_6+\\
+(p_{21}p_{62}-p_{22}p_{61}+p_{23}p_{64}-p_{24}p_{63}-p_{25}p_{66}+p_{26}p_{65})\overline{x}'_2\overline{x}'_6+\\
+(p_{31}p_{62}-p_{32}p_{61}+p_{33}p_{64}-p_{34}p_{63}-p_{35}p_{66}+p_{36}p_{65})\overline{x}'_3\overline{x}'_6+\\
+(p_{41}p_{62}-p_{42}p_{61}+p_{43}p_{64}-p_{44}p_{63}-p_{45}p_{66}+p_{46}p_{65})\overline{x}'_4\overline{x}'_6+\\
+(p_{51}p_{62}-p_{52}p_{61}+p_{53}p_{64}-p_{54}p_{63}-p_{55}p_{66}+p_{56}p_{65})\overline{x}'_5\overline{x}'_6+\\
+f(\overline{x}'_1,\overline{x}'_2,\overline{x}'_3,\overline{x}'_4,\overline{x}'_5)=0,
\end{gather*}
где $f(\overline{x}'_1,\overline{x}'_2,\overline{x}'_3,\overline{x}'_4,\overline{x}'_5)$ -- некоторый многочлен от переменных $\overline{x}'_1$, $\overline{x}'_2$, $\overline{x}'_3$, $\overline{x}'_4$, $\overline{x}'_5$, не содержащий переменную~$\overline{x}'_6$.
Так как
$\overline{x}'_1\overline{x}'_6$,
$\overline{x}'_2\overline{x}'_6$,
$\overline{x}'_3\overline{x}'_6$,
$\overline{x}'_4\overline{x}'_6$,
$\overline{x}'_5\overline{x}'_6$,
линейно независимы, то
\begin{gather*}
p_{11}p_{62}-p_{12}p_{61}+p_{13}p_{64}-p_{14}p_{63}-p_{15}p_{66}+p_{16}p_{65}=0;\\
p_{21}p_{62}-p_{22}p_{61}+p_{23}p_{64}-p_{24}p_{63}-p_{25}p_{66}+p_{26}p_{65}=0;\\
p_{31}p_{62}-p_{32}p_{61}+p_{33}p_{64}-p_{34}p_{63}-p_{35}p_{66}+p_{36}p_{65}=0;\\
p_{41}p_{62}-p_{42}p_{61}+p_{43}p_{64}-p_{44}p_{63}-p_{45}p_{66}+p_{46}p_{65}=0;\\
p_{51}p_{62}-p_{52}p_{61}+p_{53}p_{64}-p_{54}p_{63}-p_{55}p_{66}+p_{56}p_{65}=0.\\
\end{gather*}
Добавив очевидное равенство, получим
\begin{gather*}
p_{11}p_{62}-p_{12}p_{61}+p_{13}p_{64}-p_{14}p_{63}-p_{15}p_{66}+p_{16}p_{65}=0;\\
p_{21}p_{62}-p_{22}p_{61}+p_{23}p_{64}-p_{24}p_{63}-p_{25}p_{66}+p_{26}p_{65}=0;\\
p_{31}p_{62}-p_{32}p_{61}+p_{33}p_{64}-p_{34}p_{63}-p_{35}p_{66}+p_{36}p_{65}=0;\\
p_{41}p_{62}-p_{42}p_{61}+p_{43}p_{64}-p_{44}p_{63}-p_{45}p_{66}+p_{46}p_{65}=0;\\
p_{51}p_{62}-p_{52}p_{61}+p_{53}p_{64}-p_{54}p_{63}-p_{55}p_{66}+p_{56}p_{65}=0;\\
p_{61}p_{62}-p_{62}p_{61}+p_{63}p_{64}-p_{64}p_{63}-p_{65}p_{66}+p_{66}p_{65}=0.
\end{gather*}
Отсюда
$$
(p_{62}, -p_{61}, p_{64}, -p_{63}, -p_{66}, p_{65})\cdot P^{T}=0.
$$
Матрица $P$ невырожденная и, в частности, не имеет нулевых строк, а значит, полученное равенство невозможно.
Противоречие доказывает предложение. \end{proof}

Из следствия~\ref{cor1} и предложения~\ref{pr3} вытекает справедливость следующего утверждения. 

\begin{cor}\label{cor2}Если в многообразии колец~$\mathfrak{M}$ все конечные кольца однозначно определяются своими графами делителей нуля, то для любого простого числа $p$ многообразие~$\mathfrak{M}$ не содержит многообразия~$\mathfrak{M}_{1,p}$. \end{cor}

Далее, покажем, что алгебры $A_{2,p}$ и $B_{2,p}$ не изоморфны для любого простого нечетного~$p$, однако $\Gamma(A_{2,p})\cong\Gamma(B_{2,p})$. Для этого нам понадобятся некоторые вспомогательные леммы.

\begin{lemma}\label{l1}Множество $C_2=\{\overline{x}_k^{~2}, \overline{x}_i\overline{x}_j;~ (i,j)\neq(3,4), 1\leq i<j\leq 6, 1\leq k\leq6\}$ является базисом алгебры $A_{2,p}^2$ ($p>2$.)\end{lemma}
\begin{proof} Множество~$C_2$ является системой образующих векторного пространства~$A_{2,p}^2.$ Если оно линейно зависимо, то  существуют элементы $\alpha, \gamma_1, \ldots, \gamma_6, \delta_{ij}\in\mathbb{Z}_p,$ не все равные нулю, такие, что в алгебре~$F_2$ справедливо  равенство \begin{equation}\sum\limits^{6}_{i=1}\gamma_ix_i^2+\sum\limits_{\tiny\begin{array}[c]{c}  i<j \\ (i,j)\neq (3,4)\end{array}}\delta_{ij}x_ix_j=\alpha(x_3x_4-x_1x_2).\end{equation} Полагая в равенстве~(1) $x_1=\ldots=x_{i-1}=x_{i+1}=\ldots=x_6,$ где $1\leq i\leq 6,$ получим, что $\gamma_ix_i^2=0,$ т.е. $\gamma_1=\ldots=\gamma_6=0.$ Положим теперь в равенстве~(1) ${x_3=x_4=x_5=x_6=0.}$ Тогда $\delta_{12}=-\alpha,$ другими словами, $\sum\limits_{\tiny\begin{array}[c]{c}  i<j \\ (i,j)\neq (3,4)\\ (i,j)\neq (1,2) \end{array}}\delta_{ij}x_ix_j=\alpha x_3x_4.$ Положим, наконец, ${x_1=x_2=x_5=x_6=0.}$ Получим, что $\alpha=0$ и $\delta_{ij}=0$ для всех $i,j,$ таких, что $i<j $ и $(i,j)\neq (3,4), (1,2)$. \end{proof}

\begin{lemma}\label{l2}Если $a=\sum\limits^{6}_{i=1}\alpha_i\overline{x}_i+u\in A_{2,p},$ где $p>2,$ $\alpha_i\in\mathbb{Z}_p$ для всех чисел $i\in\{1,\ldots,6\},$ $u\in A_{2,p}^2,$ $a\notin A_{2,p}^2,$ то $ann(a)=A_{2,p}^2.$\end{lemma}
\begin{proof} Пусть $b=\sum\limits^{6}_{i=1}\beta_i\overline{x}_i+v\in ann(a),$ где $\beta_i\in\mathbb{Z}_p,$  $i\in\{1,\ldots,6\},$ $v\in A_{2,p}^2.$ Тогда  
\begin{equation*}\begin{split}
\overline{0}=ab=& \left(\sum\limits^{6}_{i=1}\alpha_i\overline{x}_i\right)\left(\sum\limits^{6}_{i=1}\beta_i\overline{x}_i\right)=\\=&\sum\limits^{6}_{i=1}\alpha_i\beta_i\overline{x}_i^{~2}+\left(\alpha_1\beta_2+\alpha_2\beta_1+\alpha_3\beta_4+\alpha_4\beta_3\right)\overline{x}_1\overline{x}_2+\sum\limits_{\tiny\begin{array}[c]{c}  i<j \\ (i,j)\neq (1,2)\\ (i,j)\neq (3,4)\end{array}}\left(\alpha_i\beta_j+\alpha_j\beta_i\right)\overline{x}_i\overline{x}_j.\end{split}\end{equation*} Из леммы~\ref{l1} следует, что $$\alpha_1\beta_1=\ldots=\alpha_6\beta_6=0, \alpha_1\beta_2+\alpha_2\beta_1+\alpha_3\beta_4+\alpha_4\beta_3=0 \mbox{~и~} \alpha_i\beta_j+\alpha_j\beta_i=0$$ для всех $i,j\in\{1,\ldots,6\},$ таких, что $i<j,$ $(i,j)\neq (1,2)$ и $(i,j)\neq (3,4)$. 

Пусть $\alpha_1\neq 0.$ Тогда из равенств $\alpha_1\beta_1=0$ и $\alpha_1\beta_i+\alpha_i\beta_1=0,$ где $i=3,4,5,6,$ получаем, что $\beta_1=0$ и $\beta_3=\beta_4=\beta_5=\beta_6=0.$ Далее, из равенства $\alpha_1\beta_2+\alpha_2\beta_1+\alpha_3\beta_4+\alpha_4\beta_3=0$ следует, что $\alpha_1\beta_2=0,$ т.е. $\beta_2=0.$ Другими словами, мы показали, что в этом случае $b=v\in A_{2,p}^2.$

Расссмотрим теперь случай, когда $\alpha_1=0.$ Ввиду доказанного выше, можем считать, что $\alpha_2=\alpha_3=\alpha_4=0.$ Значит, имеем равенства $\alpha_5\beta_5=\alpha_6\beta_6=0,$ $\alpha_5\beta_6+\alpha_6\beta_5=0.$ Поскольку $a\notin A_{2,p}^2,$ то $\alpha_6\neq 0$ или $\alpha_5\neq 0.$ Пусть $\alpha_5\neq 0.$ Тогда из равенства ${\alpha_5\beta_5=0}$ получаем, что $\beta_5=0.$ Поэтому $\alpha_5\beta_6=0.$ Отсюда $\beta_6=0.$ Другими словами, $b\in A_{2,p}^2.$ Случай, когда $\alpha_6\neq 0$ рассматривается аналогично. \end{proof}

Таким образом, граф $\Gamma(A_{2,p})$ имеет следующее строение: множество ненулевых элементов из множества~$A_{2,p}^2$ образуют полный подграф, а любая вершина $a\in A_{2,p}\setminus A_{2,p}^2$ смежна со всеми вершинами из этого подграфа и только с ними.

Для алгебры~$B_{2,p}$ ($p>2$) справедливы аналоги лемм~\ref{l1}--\ref{l2}.

\begin{lemma}\label{l3}Множество $D_2=\{\overline{x}_k^{~2}, \overline{x}_i\overline{x}_j;~ (i,j)\neq(5,6), 1\leq i<j\leq 6, 1\leq k\leq6\}$ является базисом алгебры $B_{2,p}^2$ ($p>2$.)\end{lemma}

\textit{Доказательство} аналогично доказательству леммы~\ref{l1}.

\begin{lemma}\label{l4} Пусть $a=\sum\limits^{6}_{i=1}\alpha_i\overline{x}_i+u\in B_{2,p},$ где $p>2,$ $\alpha_i\in\mathbb{Z}_p$ для всех чисел $i\in\{1,\ldots,6\},$ $u\in B_{2,p}^2$ и $a\notin B_{2,p}^2.$ Тогда $ann(a)=B_{2,p}^2.$\end{lemma}
\begin{proof} Возьмем элемент $b=\sum\limits^{6}_{i=1}\beta_i\overline{x}_i+v\in ann(a),$ где $\beta_i\in\mathbb{Z}_p$  для всех чисел $i\in\{1,\ldots,6\}$ и $v\in A_{2,p}^2.$ Получаем
\begin{equation*}\begin{split}
\overline{0}=ab=\sum\limits^{6}_{i=1}\alpha_i\beta_i\overline{x}_i^{~2}&+\left(\alpha_1\beta_2+\alpha_2\beta_1+\alpha_5\beta_6+\alpha_6\beta_5\right)\overline{x}_1\overline{x}_2+\\&+\left(\alpha_3\beta_4+\alpha_4\beta_3+\alpha_5\beta_6+\alpha_6\beta_5\right)\overline{x}_3\overline{x}_4+\sum\limits_{\tiny\begin{array}[c]{c}  i<j \\ (i,j)\neq (1,2)\\ (i,j)\neq (3,4)\\ (i,j)\neq (5,6)\end{array}}\left(\alpha_i\beta_j+\alpha_j\beta_i\right)\overline{x}_i\overline{x}_j.\end{split}\end{equation*} Из леммы~\ref{l3} следует, что 
\begin{equation*}\begin{split}\alpha_1\beta_1=\ldots=\alpha_6\beta_6=0,& \alpha_1\beta_2+\alpha_2\beta_1+\alpha_5\beta_6+\alpha_6\beta_5=0,\\& \alpha_3\beta_4+\alpha_4\beta_3+\alpha_5\beta_6+\alpha_6\beta_5=0 \mbox{~и~} \alpha_i\beta_j+\alpha_j\beta_i=0\end{split}\end{equation*}  для всех $i,j\in\{1,\ldots,6\},$ таких, что $i<j$ и $(i,j)\neq (1,2),(3,4),(5,6)$.

Предположим, что $\alpha_1\neq 0.$ Из равенства $\alpha_1\beta_1=0$ получаем, что $\beta_1=0.$ Далее, из равенств $\alpha_1\beta_i+\alpha_i\beta_1=0,$ $i=3,4,5,6,$ следует, что $\beta_3=\beta_4=\beta_5=\beta_6=0.$ Наконец, из равенства $\alpha_1\beta_2+\alpha_2\beta_1+\alpha_5\beta_6+\alpha_6\beta_5=0$ вытекает, что $\beta_2=0.$ Таким образом, $b=v\in B_{2,p}^2.$ 

Аналогично рассматриваются случаи, когда отличен от нуля один из коэффициентов $\alpha_2, \alpha_3, \alpha_4.$ Итак, мы можем полагать, что $\alpha_i=0$ для $i=1,2,3,4.$ Так как $a\notin A_{2,p}^2,$ то $\alpha_6\neq 0$ или $\alpha_5\neq 0.$ Предположим, что $\alpha_5\neq 0.$ Из равенства $\alpha_5\beta_5=0$ получаем, что $\beta_5=0.$ Следовательно, $\alpha_5\beta_6=0,$ т.е. и $\beta_6=0.$ Из равенств $\alpha_i\beta_5+\alpha_5\beta_i=0,$ где $i=1,2,3,4,$ следует, что $\beta_i=0$ для $i=1,2,3,4.$ Значит, $b=v\in B_{2,p}^2.$\end{proof}

Из лемм~\ref{l2} и~\ref{l4} получаем

\begin{cor}\label{cor3}$\Gamma(A_{2,p})\cong\Gamma(B_{2,p})$ для любого простого нечетного числа~$p$.\end{cor}

\begin{pr}\label{pr4}Алгебра $A_{2,p}$ не изоморфна алгебре $B_{2,p}$  для любого простого нечетного числа~$p$.\end{pr}
\begin{proof} Предположим противное. Пусть $\varphi$  -- изоморфизм $B_{2,p}$ на $A_{2,p}$. 
Обозначим
$\overline{x_i}=x_i+\langle x_1x_2+x_3x_4-x_5x_6\rangle\in B_{2,p}$ и
$\overline{x}'_i=x_i+\langle x_1x_2-x_3x_4\rangle\in A_{2,p}$, $i=\overline{1,6}$.
Заметим что $\varphi(B_{2,p}^2)\subseteq A_{2,p}^2$, а значит, отображение $\overline{\varphi}: B_{2,p}/B_{2,p}^2 \rightarrow A_{2,p}/A_{2,p}^2$, определенное по правилу
$\overline{\varphi}(a+B_{2,p}^2)=\varphi(a)+A_{2,p}^2$, $a\in B_{2,p}$, также является изоморфизмом. 
Отсюда следует, что существует невырожденная матрица $P=(p_{ij})_{6\times 6}$, $p_{ij}\in \mathbb{Z}_p,$ такая, что 
$$
\left\{
  \begin{array}{ll}
\varphi(\overline{x}_1)=p_{11}\overline{x}'_1+p_{21}\overline{x}'_2+p_{31}\overline{x}'_3+p_{41}\overline{x}'_4+p_{51}\overline{x}'_5+p_{61}\overline{x}'_6+b_1;\\
\varphi(\overline{x}_2)=p_{12}\overline{x}'_1+p_{22}\overline{x}'_2+p_{32}\overline{x}'_3+p_{42}\overline{x}'_4+p_{52}\overline{x}'_5+p_{62}\overline{x}'_6+b_2;\\
\varphi(\overline{x}_3)=p_{13}\overline{x}'_1+p_{23}\overline{x}'_2+p_{33}\overline{x}'_3+p_{43}\overline{x}'_4+p_{53}\overline{x}'_5+p_{63}\overline{x}'_6+b_3;\\
\varphi(\overline{x}_4)=p_{14}\overline{x}'_1+p_{24}\overline{x}'_2+p_{34}\overline{x}'_3+p_{44}\overline{x}'_4+p_{54}\overline{x}'_5+p_{64}\overline{x}'_6+b_4;\\
\varphi(\overline{x}_5)=p_{15}\overline{x}'_1+p_{25}\overline{x}'_2+p_{35}\overline{x}'_3+p_{45}\overline{x}'_4+p_{55}\overline{x}'_5+p_{65}\overline{x}'_6+b_5;\\
\varphi(\overline{x}_6)=p_{16}\overline{x}'_1+p_{26}\overline{x}'_2+p_{36}\overline{x}'_3+p_{46}\overline{x}'_4+p_{56}\overline{x}'_5+p_{66}\overline{x}'_6+b_6,
  \end{array}
\right.
$$
где $b_1$, $b_2$, $b_3$, $b_4$, $b_5$, $b_6$ -- некоторые элементы из $A_{2,p}^2$. 

Учитывая свойства изоморфизма, имеем
\begin{gather*}
\varphi(\overline{x}_1\overline{x}_2+\overline{x}_3\overline{x}_4-\overline{x}_5\overline{x}_6)=\\
=(p_{11}\overline{x}'_1+p_{21}\overline{x}'_2+p_{31}\overline{x}'_3+p_{41}\overline{x}'_4+p_{51}\overline{x}'_5+p_{61}\overline{x}'_6)
(p_{12}\overline{x}'_1+p_{22}\overline{x}'_2+p_{32}\overline{x}'_3+p_{42}\overline{x}'_4+p_{52}\overline{x}'_5+p_{62}\overline{x}'_6)+\\
+(p_{13}\overline{x}'_1+p_{23}\overline{x}'_2+p_{33}\overline{x}'_3+p_{43}\overline{x}'_4+p_{53}\overline{x}'_5+p_{63}\overline{x}'_6)
(p_{14}\overline{x}'_1+p_{24}\overline{x}'_2+p_{34}\overline{x}'_3+p_{44}\overline{x}'_4+p_{54}\overline{x}'_5+p_{64}\overline{x}'_6)-\\
-(p_{15}\overline{x}'_1+p_{25}\overline{x}'_2+p_{35}\overline{x}'_3+p_{45}\overline{x}'_4+p_{55}\overline{x}'_5+p_{65}\overline{x}'_6)
(p_{16}\overline{x}'_1+p_{26}\overline{x}'_2+p_{36}\overline{x}'_3+p_{46}\overline{x}'_4+p_{56}\overline{x}'_5+p_{66}\overline{x}'_6)=\\
=(p_{11}p_{62}+p_{12}p_{61}+p_{13}p_{64}+p_{14}p_{63}-p_{15}p_{66}-p_{16}p_{65})\overline{x}'_1\overline{x}'_6+\\
+(p_{21}p_{62}+p_{22}p_{61}+p_{23}p_{64}+p_{24}p_{63}-p_{25}p_{66}-p_{26}p_{65})\overline{x}'_2\overline{x}'_6+\\
+(p_{31}p_{62}+p_{32}p_{61}+p_{33}p_{64}+p_{34}p_{63}-p_{35}p_{66}-p_{36}p_{65})\overline{x}'_3\overline{x}'_6+\\
+(p_{41}p_{62}+p_{42}p_{61}+p_{43}p_{64}+p_{44}p_{63}-p_{45}p_{66}-p_{46}p_{65})\overline{x}'_4\overline{x}'_6+\\
+(p_{51}p_{62}+p_{52}p_{61}+p_{53}p_{64}+p_{54}p_{63}-p_{55}p_{66}-p_{56}p_{65})\overline{x}'_5\overline{x}'_6+\\
+(p_{61}p_{62}+p_{63}p_{64}-p_{65}p_{66})\overline{x}'_6\overline{x}'_6+f(\overline{x}'_1,\overline{x}'_2,\overline{x}'_3,\overline{x}'_4,\overline{x}'_5)=0,
\end{gather*}
где $f(\overline{x}'_1,\overline{x}'_2,\overline{x}'_3,\overline{x}'_4,\overline{x}'_5)$ -- некоторый многочлен от переменных $\overline{x}'_1$, $\overline{x}'_2$, $\overline{x}'_3$, $\overline{x}'_4$, $\overline{x}'_5$, не содержащий $\overline{x}'_6$.
Поскольку элементы
$\overline{x}'_1\overline{x}'_6$,
$\overline{x}'_2\overline{x}'_6$,
$\overline{x}'_3\overline{x}'_6$,
$\overline{x}'_4\overline{x}'_6$,
$\overline{x}'_5\overline{x}'_6$,
$\overline{x}'_6\overline{x}'_6$ линейно независимы, то
\begin{gather*}
p_{11}p_{62}+p_{12}p_{61}+p_{13}p_{64}+p_{14}p_{63}-p_{15}p_{66}-p_{16}p_{65}=0;\\
p_{21}p_{62}+p_{22}p_{61}+p_{23}p_{64}+p_{24}p_{63}-p_{25}p_{66}-p_{26}p_{65}=0;\\
p_{31}p_{62}+p_{32}p_{61}+p_{33}p_{64}+p_{34}p_{63}-p_{35}p_{66}-p_{36}p_{65}=0;\\
p_{41}p_{62}+p_{42}p_{61}+p_{43}p_{64}+p_{44}p_{63}-p_{45}p_{66}-p_{46}p_{65}=0;\\
p_{51}p_{62}+p_{52}p_{61}+p_{53}p_{64}+p_{54}p_{63}-p_{55}p_{66}-p_{56}p_{65}=0;\\
p_{61}p_{62}+p_{63}p_{64}-p_{65}p_{66}=0.
\end{gather*}
Изменив последнее уравнение, получим
\begin{gather*}
p_{11}p_{62}+p_{12}p_{61}+p_{13}p_{64}+p_{14}p_{63}-p_{15}p_{66}-p_{16}p_{65}=0;\\
p_{21}p_{62}+p_{22}p_{61}+p_{23}p_{64}+p_{24}p_{63}-p_{25}p_{66}-p_{26}p_{65}=0;\\
p_{31}p_{62}+p_{32}p_{61}+p_{33}p_{64}+p_{34}p_{63}-p_{35}p_{66}-p_{36}p_{65}=0;\\
p_{41}p_{62}+p_{42}p_{61}+p_{43}p_{64}+p_{44}p_{63}-p_{45}p_{66}-p_{46}p_{65}=0;\\
p_{51}p_{62}+p_{52}p_{61}+p_{53}p_{64}+p_{54}p_{63}-p_{55}p_{66}-p_{56}p_{65}=0;\\
p_{61}p_{62}+p_{62}p_{61}+p_{63}p_{64}+p_{64}p_{63}-p_{65}p_{66}-p_{66}p_{65}=0.
\end{gather*}
Отсюда
$$
(p_{62}, p_{61}, p_{64}, p_{63}, -p_{66}, -p_{65})\cdot P^{T}=0.
$$
Так как матрица $P$ невырожденная, то полученное равенство невозможно.
Противоречие. \end{proof}

Из следствия~\ref{cor3} и предложения~\ref{pr4} вытекает справедливость следующего утверждения.

\begin{cor}\label{cor4}Если в многообразии колец~$\mathfrak{M}$ все конечные кольца однозначно определяются своими графами делителей нуля, то для любого нечетного простого числа $p$ многообразие~$\mathfrak{M}$ не содержит многообразия~$\mathfrak{M}_{2,p}$. \end{cor}

Пусть $\mathfrak{M}$~--- многообразие, в котором все конечные кольца однозначно определяются своими графами делителей нуля. В работе~\cite{semr} доказано, что $T(\mathfrak{M})$ содержит многочлены вида $mx,$ ${d x+ x^2 g(x),}$ причем $m=p_1^{\alpha_1}\ldots p_s^{\alpha_s},$ $\alpha_i\leq 3$ для всех $i\leq s,$  $g(x)\in \mathbb{Z}[x]$, $d=1$ или ${d=p_{i_1}\ldots p_{i_k},}$ где $p_{i_1},\ldots, p_{i_k}$~-- попарно различные простые делители числа~$m$. 

Далее, пусть $\mathfrak{N}_{i}=var\left\langle p_i^{\alpha_i}x=0\right\rangle\cap \mathfrak{M},$ где $1\leq i\leq s.$ Тогда $\mathfrak{M}=\mathfrak{N}_{1}\vee \ldots \vee\mathfrak{N}_{s}$ и $T(\mathfrak{M})=T(\mathfrak{N}_{1})\cap\ldots\cap T(\mathfrak{N}_{s}).$

\begin{lemma}\label{l5} Если кольцо $A$ удовлетворяет тождеству $x(y-y^n)=0,$ $n\geq 2,$ а кольцо~$B$ удовлетворяет тождеству $x(y-y^m)=0,$ $m\geq 2,$ то в кольце $A\oplus B$  выполняется тождество $$x(y-y^{(n-1)(m-1)+1})=0.$$\end{lemma}
\begin{proof} Пусть $a,b$~--- произвольные элементы из кольца~$A$. Тогда $$ab=ab^n=ab\cdot b^{n-1}= (ab\cdot b^{n-1})b^{n-1}=ab\cdot b^{2(n-1)}=\ldots=ab^{1+(n-1)t},$$ где $t$~--- произвольное неотрицательное целое число. Другими словами, в кольце~$A$ выполняются тождества $xy=xy^{1+(n-1)t},$  $t\geq 1.$ Аналогично доказывается, что кольцо $B$ удовлетворяет тождествам $xy=xy^{1+(m-1)s},$ где $s$~--- произвольное целое неотрицательное число. Значит, $xy=xy^{1+(n-1)(m-1)}$~--- тождество в кольце~$A\oplus B$. \end{proof}

\vspace{2cm}

\begin{pr}\label{pr5} $T(\mathfrak{N}_i)$ содержит многочлен $x(y-y^{p_i}),$ $i=1,\ldots, s.$ \end{pr}

\begin{proof}Из работы~\cite{aejm2012} следует, что $T(\mathfrak{N}_i)$ содержит многочлены вида $p_i^{\alpha_i}x$ и $dx+x^2f(x),$ где $d=1$ или $d=p_i,$ $f(x)\in\mathbb{Z}[x],$ $\alpha_i\leq3.$ Если $d=1,$ то идеал тождеств $T(\mathfrak{N}_i)$ содержит многочлен  $x+x^2f(x)$, и из~\cite{Jacobson} следует, что $\mathfrak{N}_i$ порождается полем $\mathbb{Z}_{p_i}$ или $\mathfrak{N}_i=var\left\langle x=0\right\rangle.$ В каждом из этих случае получаем, что $T(\mathfrak{N}_i)$ содержит многочлен $x(y-y^{p_i}).$ 

Рассмотрим случай, когда $d=p_i,$ т.е. $p_ix+x^2f(x)\in T(\mathfrak{N}_i).$ Из следствия~\ref{cor2} имеем, что многообразие~$\mathfrak{N}_i$ не содержит многообразия~$\mathfrak{M}_{1,p_i}$. Это означает, что существует многочлен $g(x_1, \ldots, x_N)\in T(\mathfrak{N}_i), $ существенно зависящий от переменных $x_1, \ldots, x_N,$ такой, что $g(x_1, \ldots, x_N)\notin \{xyz, x^2, p_ix\}^T.$ Ясно, что $N\leq 2.$ 

Пусть $N=1$. В этом случае многочлен $g$ можно записать в виде: $$g(x)=bx+ap_ix+x^2h(x),$$ где $h(x)\in\mathbb{Z}[x],$ $a,b\in \mathbb{Z},$ причем число~$b$ не делится на $p_i$,  По лемме о НОД существуют целые числа~$u,v$, такие, что $bu+p_iv=1.$ Поскольку $p_ix+x^2f(x)\in T(\mathfrak{N}_i),$ то $T(\mathfrak{N}_i)$ содержит многочлен вида $g_1(x)=x+cp_ix+x^2h_1(x)$ для некоторых $c\in\mathbb{Z}$ и ${h_1(x)\in\mathbb{Z}[x].}$ Отсюда, наконец, получаем, что $$\left(x+cp_ix+x^2h_1(x)\right)-c\left(p_ix+x^2f(x)\right)=x+x^2(h_1(x)-cf(x))\in\mathfrak{N_i}.$$ Из~\cite{aejm2012} следует, что $T(\mathfrak{N}_i)$ содержит многочлен $x(y-y^{p_i}).$

Пусть теперь $N=2.$ Тогда можем записать $$g(x,y)=\alpha xy+\beta (x\circ y)+\gamma p_ixy+\varphi(x,y),$$ где $\varphi(x,y)\in \mathbb{Z}[x,y],$ $\alpha,\beta,\gamma\in\mathbb{Z},$ причем число~$\alpha$ не делится на~$p_i$, а нижняя степень многочлена $\varphi(x,y)$ больше~$2$. По лемме о НОД существуют целые числа~$u,v$, такие, что $\alpha u+p_iv=1.$ Поэтому $T(\mathfrak{N}_i)$ содержит многочлен 
\begin{equation*}\begin{split}g_1(x,y)=&ug(x,y)= u\alpha xy+u\beta (x\circ y)+u\gamma p_ixy+u\varphi(x,y)=\\=&xy+u\beta (x\circ y)+(u\gamma -v)p_ixy+u\varphi(x,y)=\\=&xy+u\beta (x\circ y)+(u\gamma -v)\left(-x^2y^2f(xy)\right)+u\varphi(x,y)=\\=& xy+\beta_1 (x\circ y)+\varphi_1(x,y),\end{split}\end{equation*} где $\beta_1=u \beta$ и $\varphi_1(x,y)$~--- некоторый многочлен с целыми коэффициентами, нижняя степень которого больше~$2$. Итак,  $T(\mathfrak{N}_i)$ содержит многочлен $$g_1(x,y)=xy+\beta_1 (x\circ y)+\varphi_1(x,y),$$ где $\beta_1\in \mathbb{Z}$  и $\varphi_1(x,y)$~--- многочлен с целыми коэффициентами, нижняя степень которого больше~$2$.

Рассмотрим случай, когда $p_i=2.$ Тогда $$g_1(x,x)=x^2+2\beta_1x^2+\varphi_1(x,x)=x^2+\beta_1\left(-x^4f(x^2)\right)+\varphi_1(x,x)=x^2+x^3 \varphi_2(x)\in T(\mathfrak{N}_{i})$$ для некоторого многочлена $\varphi_2(x)\in\mathbb{Z}[x].$ Отсюда следует, что многообразие~$\mathfrak{N}_i$ удовлетворяет тождеству вида $x\circ y+\psi(x,y)=0$, где $\psi(x,y)\in\mathbb{Z}[x,y],$ причем нижняя степень многочлена~$\psi(x,y)$ больше~$2$. Следовательно, $$g_1(x,y)=xy+\beta_1\left(-\psi(x,y)\right)+\varphi_1(x,y)=xy+\psi_1(x,y),$$ где $\psi_1(x,y)=-\beta_1\psi(x,y)+\varphi_1(x,y)$~--- многочлен, нижняя степень которого больше~$2$. Из теоремы~1 работы~\cite{aejm2012} имеем, что $\mathfrak{N}_i\subseteq var~N_{0,2}\oplus \mathbb{Z}_2.$ Следовательно, многообразие~$\mathfrak{N}_i$ удовлетворяет тождеству ~$x(y-y^{2})=0$.

Рассмотрим теперь случай, когда $p_i$~--- нечетное число. Из следствия~\ref{cor4} имеем, что многообразие~$\mathfrak{N}_i$ не содержит многообразия~$\mathfrak{M}_{2,p_i}$. Значит, существует многочлен $p(x_1, \ldots, x_M)\in T(\mathfrak{N}_i), $ существенно зависящий от переменных $x_1, \ldots, x_M,$ такой, что $p(x_1, \ldots, x_M)\notin \{xyz, [x,y], p_ix\}^T.$ Ясно, что $M\leq 2.$ 

Положим сначала $M=1.$ Тогда можем записать $$p(x)=\lambda x+\mu p_ix +\nu x^2+\delta p_i x^2+x^3\sigma(x),$$ где $\sigma(x)\in\mathbb{Z}[x],$ $\lambda,\mu,\nu,\delta \in\mathbb{Z},$ причем либо $\lambda$~--- ненулевое число, не делящееся на~$p_i$, либо $\nu$~--- ненулевое число, не делящееся на~$p_i$. Пусть $\lambda$ не равно нулю и взаимно просто с числом~$p_i$. Тогда лемме о НОД существуют целые числа $u,v$, такие, что $\lambda u+ p_i v=1.$ Отсюда следует, что 
\begin{equation*}\begin{split}up(x)=&x+p_i(\mu u- v)x +\nu ux^2+\delta up_i x^2+ux^3\sigma(x)=\\=&x+(\mu u- v)(-x^2f(x)) +\nu ux^2+\delta up_i x^2+ux^3\sigma(x),\end{split}\end{equation*} т.е. $\mathfrak{N}_i$ удовлетворяет тождеству вида $x+x^2q(x)=0,$ где $q(x)\in\mathbb{Z}[x].$ Ранее мы отмечали, что в этом случае $x(y-y^{p_i})\in T(\mathfrak{N}_i).$ Пусть теперь $\lambda=0.$ Тогда число~$\nu$ не равно нулю и не делится на~$p_i$. Отсюда следует, что лемме о НОД существуют целые числа $u,v$, такие, что $\nu u+ p_i v=1.$  Значит, имеем
\begin{equation*}\begin{split}up(x)=&\mu up_ix +x^2+p_i(\delta u -v)x^2+ux^3\sigma(x)= \\=&\mu up_ix +x^2+(\delta u -v)(-x^4f(x^2))+ux^3\sigma(x)= \mu up_ix +x^2+x^3\sigma_1(x)\end{split}\end{equation*}  для некоторого многочлена $\sigma_1(x)\in\mathbb{Z}[x].$ Линеаризуя многочлен $\mu up_ix +x^2+x^3\sigma_1(x),$ мы получим, что идеал тождеств~$T(\mathfrak{N}_i)$ содержит многочлен вида $x\circ y+w(x,y),$ где $w(x,y)$~--- некоторый многочлен с целыми коэффициентами, нижняя степень которого больше~$2$. Из тождеств $x\circ y+w(x,y)=0$ и $g_1(x,y)=0$ получаем, что $T(\mathfrak{N}_i)$ содержит многочлен вида~$xy+w_1(x,y),$ где $w_1(x,y)\in\mathbb{Z}[x,y]$ и нижняя степень многочлена~$w_1(x,y)$ больше~$2$. По теореме~1 статьи~\cite{aejm2012} следует, что многообразие~$\mathfrak{N}_i$ удовлетворяет тождеству $x(y-y^{p_i})=0.$

Пусть, наконец, $M=2.$ Тогда можем записать $$p(x,y)=\lambda_1xy+\mu_1[x,y]+\nu_1p_ixy+\omega(x,y),$$ где $\lambda_1, \mu_1, \nu_1\in\mathbb{Z},$ причем $\lambda_1$ не делится на~$p_i$, и $\omega(x,y)$~--- многочлен с целыми коэффициентами, нижняя степень которого больше~$2$. Используя лемму о НОД так же, как это было сделано выше, мы можем считать, что $\lambda_1=1,$ т.е. $$p(x,y)=xy+\mu_1[x,y]+\nu_1p_ixy+\omega(x,y).$$ Пользуясь тождеством $p_ix+x^2f(x)=0,$ получаем, что $$p(x,y)=xy+\mu_1[x,y]+\nu_1(-x^2y^2f(xy))+\omega(x,y)=xy+\mu_1[x,y]+\omega_1(x,y)$$ для некоторого многочлена~$\omega_1(x,y)\in\mathbb{Z}[x,y],$ нижняя степень которого больше~$2$. Тогда $p(x,x)=x^2+x^3\omega_2(x),$ где $\omega_2\in\mathbb{Z}[x].$ Линеаризуя тождество $p(x,x)=0,$ получим, что многообразие~$\mathfrak{N}_i$ удовлетворяет тождеству вида $x\circ y + w(x,y)=0,$ где $w(x,y)$~--- некоторый многочлен с целыми коэффициентами, нижняя степень которого больше~$2$. Из тождеств $x\circ y + w(x,y)=0$ и $g_1(x,y)=0$ следует тождество вида $xy+w_1(x,y)=0,$ где $w_1(x,y)\in\mathbb{Z}[x,y]$ и нижняя степень многочлена~$w_1(x,y)$ больше~$2$. По теореме~1 статьи~\cite{aejm2012} получаем, что многообразие~$\mathfrak{N}_i$ удовлетворяет тождеству $x(y-y^{p_i})=0.$ \end{proof} 

Из предложения~\ref{pr5} и леммы~\ref{l5} получаем следующее утверждение.

\begin{cor}\label{cor5} Пусть $\mathfrak{M}$~--- многообразие, в котором все конечные кольца однозначно определяются своими графами делителей нуля. Тогда $T(\mathfrak{M})$ содержит многочлен вида $x(y-y^N),$ где $N\geq 2.$  \end{cor}

Теперь мы можем доказать основной результат настоящей работы.

\begin{theo}\label{th} Для любого многообразия~$\mathfrak{M}$ ассоциативных колец следующие условия эквивалентны:
\begin{enumerate}
	\item Произвольное конечное кольцо из~$\mathfrak{M}$  однозначно определяется своим графом делителей нуля;
		\item $\mathfrak{M}\subseteq var\left\langle N_{0,p_1}\oplus \ldots \oplus N_{0, p_s}\oplus \mathbb{Z}_p\right\rangle,$ где $s\in \mathbb{N}$ и $(p_i,p)\neq (3,2)$ для любого числа $i\in\{1,\ldots, s\}.$
\end{enumerate}
\end{theo}
\begin{proof} Импликация $(2)\Rightarrow (1)$ следует из предложения~4 работы~\cite{semr}.  Докажем импликацию  $(1)\Rightarrow (2)$. Пусть в многообразии~$\mathfrak{M}$ все конечные кольца однозначно определяются своими графами делителей нуля. Тогда по следствию~\ref{cor5} идеал тождеств $T(\mathfrak{M})$ содержит многочлен вида $x(y-y^N),$ где $N\geq 2.$ По теореме~1.1 из статьи~\cite{aejm2012} имеем, что $\mathfrak{M}\subseteq var\left\langle N_{0,p_1}\oplus \ldots \oplus N_{0, p_s}\oplus \mathbb{Z}_p\right\rangle,$ где $s\in \mathbb{N}$ и $(p_i,p)\neq (3,2)$ для любого числа $i\in\{1,\ldots, s\}.$ \end{proof}

\end{document}